\documentclass[12pt]{article}

\usepackage{amsfonts,amsmath,latexsym}
\usepackage{amstext,amssymb}
\usepackage{amsthm}
\usepackage{color}
\usepackage{soul}

\usepackage{tikz}

\newcommand {\R}{\mathbb{R}}

\newcommand {\B}{\mathbf{B}}

\newtheorem {thm} {Theorem}

\newtheorem {prop} [thm] {Proposition}
\newtheorem {lemma}{Lemma}
\newtheorem {cor} [thm] {Corollary}

\newcommand{\beq}{\begin{equation}}
\newcommand{\eeq}{\end{equation}}
\numberwithin{equation}{section}
\newcommand{\cpr}{\mathcal{C}_{p,r}}

\begin {document}

\title {On the rate of change of the best constant in the Sobolev inequality} 
\author{Tom Carroll\footnote{School of Mathematical Sciences, 
University College Cork, Ireland, {\tt t.carroll@ucc.ie}},
Mouhamed Moustapha Fall\footnote{African Institute for Mathematical Sciences 
of Senegal, {\tt mouhamed.m.fall@aims-senegal.org}}, and 
Jesse Ratzkin\footnote{Department of Mathematics and Applied Mathematics, University of 
Cape Town, South Africa, {\tt jesse.ratzkin@uct.ac.za}} } 


\maketitle 

\begin {abstract}\noindent We estimate the rate of change of the 
best constant in the Sobolev inequality of a Euclidean domain which moves outward. 
Along the way we prove an inequality which reverses the usual H\"older inequality, 
which may be of independent interest. 
\end {abstract} 

\section {Introduction} 
The Sobolev inequality, in its many and varied forms,  is a key functional 
geometric inequality by which integrability properties of a function are inferred 
from integrability properties of its derivative. 
In $n$ dimensions, $n \geq 2$, and for $r\in [1, n)$, the most basic form 
of the inequality states that there 
is a finite constant  $S_r(\R^n)$ such that for any real-valued smooth function 
of compact support in $\R^n$, 
\begin{equation}
\label{tomsobolev}
\Vert u \Vert_{L^{r^*}(\R^n)}  \leq S_r(\R^n) \Vert \nabla u \Vert _{L^r(\R^n)}, 
\quad r^* = \frac{nr}{n-r}.
\end{equation}
Inequalities of this form having been obtained in various settings.
It is subsequently of relevance to determine, if possible, 
the best constants in the inequalities as well as the extremal functions. 
For example, the case $r=1$ of \eqref{tomsobolev} is equivalent
to the isoperimetric inequality, and the best constant in \eqref{tomsobolev} when $r=1$ 
\emph{is\/} the isoperimetric constant -- this fact is due independently to Federer and Fleming 
and to Maz'ya, as described by Chavel \cite{ChavelIso}. 

In the setting of any open region $\Omega$ of finite volume in $\R^n$, 
it is  a consequence of the basic Sobolev inequality \eqref{tomsobolev} 
that, for $r \in [1,n)$ and $p \in [1,r^*]$,
there is a finite constant $S_{p,r}(\Omega)$ such that
\begin{equation}
\label{tomsobolevomega}
\Vert u \Vert_{L^p(\Omega)} \leq S_{p,r}(\Omega)\Vert \nabla u \Vert _{L^r(\Omega)}
\end{equation}
for any function in the Sobolev space $W_0^{1,r}(\Omega)$. We remark that, by scale 
invariance, $S_{r^*, r} (\Omega) = S_r(\R^n)$ for any open set in $\R^n$. 
The inclusion of the Sobolev space $W^{1,r}(\R^n) \subseteq L^{r^*}(\R^n)$ 
in \eqref{tomsobolev} is not a compact embedding whereas the embedding in 
\eqref{tomsobolevomega} is compact (Rellich compactness) if $p< r^*$.
The best constant in the Sobolev inequality \eqref{tomsobolevomega}, 
now in the context of the  region $\Omega$, is in essence the number 
\begin {equation} \label{sobolev-defn} 
\mathcal{C}_{p,r} (\Omega) = \inf \left \{ \frac{\int_\Omega |\nabla u |^r \,d\mu}
{\left ( \int_\Omega |u|^p \,d\mu \right )^{r/p}} : u \in \mathcal{C}^\infty_0(\Omega), 
u \not \equiv 0 \right \},
\end {equation} 
where $d\mu$ stands for Lebesgue measure in $\R^n$. 
The reason for writing the best constant in the Sobolev inequality in this form 
is historical: in two dimensions, $\mathcal{C}_{2,2}(\Omega)$ is the 
classical Rayleigh quotient for the principal frequency or bass note of the 
planar region $\Omega$ while 
$4/\mathcal{C}_{1,2}$ corresponds to the torsional rigidity of the region, 
both important physical
concepts in the context of solid mechanics. 
P\'olya and Szeg\"o's monograph \cite{PS} is a standard reference from this viewpoint. 
The relationship between the best constant $S_{p,r}(\Omega)$ in the Sobolev inequality 
\eqref{tomsobolevomega} and the eigenvalue $\cpr(\Omega)$ given by \eqref{sobolev-defn} 
is then 
\[
S_{p,r}(\Omega) = \cpr(\Omega)^{-1/r}.
\]
The Sobolev inequality implies that $\cpr(\Omega)$ is positive, 
and Rellich compactness gives the existence of a nontrivial minimizer $\phi$. 
This minimizer depends on the particular region $\Omega$ in $\R^n$ 
and on the exponents $r$ and $p$. Choose  
\begin {equation} \label {sobolev-normalization} 
\phi>0 \textrm{ in } \Omega \quad\mbox{and} \quad \int_\Omega \phi^p \,d\mu = 1,
\end {equation} 
a normalization that uniquely determines $\phi$. 
The minimizer $\phi$ satisfies an  Euler-Lagrange partial differential equation 
with zero boundary values, namely 
\begin {equation} \label {sobolev-pde} 
0 = \Delta_r \phi + \mathcal{C}_{p,r} (\Omega) \phi^{p-1} = 
\operatorname{div} (|\nabla \phi|^{r-2} \nabla \phi) + \mathcal{C}_{p,r} (\Omega) 
\phi^{p-1}, \quad \left. \phi \right |_{\partial \Omega} = 0.\end {equation} 
It will also be useful to record the scaling law
\begin {equation} \label{sobolev-scaling} 
\mathcal{C}_{p,r} (R \Omega) = R^{n-r - \frac{rn}{p}} \mathcal{C}_{p,r} (\Omega), 
\end {equation} 
which is straightforward using the change of variables $y = x/R$. 

It is clear from its definition \eqref{sobolev-defn} that if 
$\tilde\Omega \subseteq \Omega$ then $\cpr(\tilde\Omega)\geq \cpr(\Omega)$,
so that bigger regions have smaller eigenvalues just as bigger drums have lower bass notes. 
Our intention herein is to quantify the rate of decrease of the eigenvalue 
$\cpr$ as the region $\Omega$ expands. 
Assuming that $\Omega$ has $\mathcal{C}^1$ boundary that 
moves with velocity $e^w \eta$, where $\eta = \eta(\zeta)$ is the unit outward 
normal to $\Omega$ 
and $w = w(\zeta)$ is a bounded, continuous function on the boundary of $\Omega$,
we denote resulting region at time $t$ by $\Omega_t$. 
Thus, $\mathcal{C}_{p,r}(\Omega_t)$ is a non-increasing function of $t$, which begs the 
following question: can one bound 
$\dot{ \mathcal{C}}_{p,r} = \left. \frac{d}{dt} \mathcal{C}_{p,r}(\Omega_t)\right |_{t=0}$? 
Below we provide some answers, particularly in the case $p=r$ and in the case $n=r=2$. 

\begin {thm} \label{thm:p=r} Let $1< p< n$. 
There is a positive constant $K$ depending only on $n$ and $p$ such that 
\begin {equation} \label{ineq:p=r}
-\left. \frac{d}{dt}  (\mathcal{C}_{p,p}(\Omega_t))^{\frac{n-p}{p(p-1)}} \right|_{t=0} 
\geq \frac{\left ( \frac{n-p}{p} \right ) K^{\frac{1}{p-1}}} 
{ \left ( \int_{\partial \Omega} e^{(1-p)w}\,d\sigma \right )^{p-1}} 
\end {equation} 
and equality can only occur if $\Omega$ is a round ball and $w$ is constant. 
Also, 
\begin {equation} \label {ineq:p=r=n} 
-\left . \frac{d}{dt} \log \mathcal{C}_{n,n} (\Omega_t)\right|_{t=0}  \geq 
\frac{ (n-1) K^{\frac{1}{n-1}}}{\left ( \int_{\partial \Omega} e^{(1-n) w}
\,d\sigma \right )^{n-1}} 
\end {equation} 
and (as before) equality can only occur if $\Omega$ is a round ball and $w$ 
is constant. 
\end {thm} 

\begin {thm} \label{thm:n=r=2}
In dimension $n=2$ we have 
\begin {equation} \label{ineq:n=r=2}
- \left. \frac{d}{dt} \log \mathcal{C}_{p,2}(\Omega_t) \right |_{t=0} 
\geq \frac{8\pi}{p} \frac{1}{\int_{\partial \Omega} e^{-w} d\sigma}
\end {equation} 
for all $p \geq 1$, and equality implies that $\Omega$ is a round disk and $w$ is constant. 
\end {thm} 

Following the work of the first and last authors in \cite{CR2}, we also obtain an 
inequality comparing 
the eigenvalue $\mathcal{C}_{p,p}$ before and after a conformal diffeomorphism. 
Here $\B$ stands for the unit ball in $\R^n$ and $\B_t$ for the ball of radius $t$.
\begin {thm} \label{thm:conformal} 
Let $n \geq 3$ and let 
$F:\B \rightarrow \R^n$ be a conformal diffeomorphism and suppose that
\[
\int_{\partial \B_t} |DF|^{n-2} \,d\sigma \geq |\partial \B_t|^{(p-1)^2}
\]
for $0<t<1$. If $1<p<n$ then  
\begin {equation} \label{ineq:conf-p=r}
\frac{d}{dt} \left [ (\mathcal{C}_{p,p} (F(\B_t)))^{\frac{n-p}{p(p-1)}} - 
(\mathcal{C}_{p,p} (\B_t))^{\frac{n-p}{p(p-1)}} \right ]\leq 0 ,
\end {equation} 
while 
\begin {equation} \label{ineq:conf-p=r=n} 
\frac{d}{dt} \log \left ( \frac{ \mathcal{C}_{n,n} (F(\B_t))}
{\mathcal{C}_{n,n} (\B_t)} \right ) \leq 0.
\end {equation} 
Equality in either case can only occur if $F(\B_t)$ is a round ball. 
\end {thm} 

These results, in the special case $p=2$, were first obtained by two of the present 
authors, Carroll and Ratzkin, \cite[Theorem~11]{CR2}.  
It is straightforward to adapt the proofs below from the Euclidean space setting  
to the setting of a general class of Riemannian manifolds in which an 
isoperimetric inequality holds. 
The discussion in \cite{CR2} provides details of this particular extension of the results. 
Also, as described therein, Theorems~\ref{thm:p=r} and \ref{thm:n=r=2} 
apply to a large collection of geometric flows, such as 
curvature flow, under appropriate convexity hypotheses. Additionally, 
one can always apply both results to Hele-Shaw flow. As described in 
\cite{Gus}, Hele-Shaw flow models a viscous fluid injected into the space between two 
plates, and  $\partial \Omega$ moves with velocity $-\nabla G$, where 
$G$ is the Green's function for the Laplacian with 
a pole inside $\Omega \subset \R^2$ corresponding to the injection site. Thus \eqref{ineq:p=r} reads 
$$- \left. \frac{d}{dt} (\mathcal{C}_{p,p} (\Omega_t))^{\frac{2-p}{p(p-1)}} \right |_{t=0}
 \geq \frac{ \left ( \frac{2-p}{p} \right ) K^{\frac{1}{p-1}}} {\left ( \int_{\partial \Omega}
 |\nabla G|^{1-p} d\sigma \right )^{p-1}}, \qquad 1<p<2,$$
 and \eqref{ineq:n=r=2} reads 
 $$-\left. \frac{d}{dt} \log (\mathcal{C}_{p,2} (\Omega_t) \right |_{t=0}  \geq 
 \frac{8\pi}{p} \frac{1}{\int_{\partial \Omega} |\nabla G|^{-1} d\sigma} , \qquad p \geq 1.$$
Theorem \ref{thm:conformal} can be viewed as a variation on the classical 
Schwarz Lemma from complex analysis. 
One might also envisage versions of the Schwarz Lemma for $n=r=2$, using 
Theorem \ref{thm:n=r=2}, but this is already done in \cite{CR1} using a different 
technique. 

One may reasonably ask what the appropriate version of Theorem~\ref{thm:p=r} might 
be when $r \rightarrow 1^+$. 
In the limit the infimum which defines the eigenvalue $\cpr(\Omega)$ 
by \eqref{sobolev-defn} is usually not attained in the Sobolev space 
$W^{1,1}_0(\Omega)$, but rather in the space of functions with bounded mean 
oscillation, and so the Hadamard variation formula \eqref{hadamard-var2} that we use 
is not valid in the case $r=1$. 
The article \cite{KF} details this phenomenon, and describes some interesting 
relations with the Cheeger constant.

Our proofs contain two ingredients: a Hadamard variation formula \eqref{hadamard-var2}, 
and an inequality \eqref{ineq:chiti2}  which reverses the usual H\"older inequality 
in the case of extremal Sobolev functions. 
We will prove a general Hadamard variation formula which is valid in all possible cases, and 
also a reverse-H\"older inequality in the case $p=r$. One can find the requisite 
reverse-H\"older inequality for the case $n=r=2$ in \cite{CR1}. 
It now seems clear that a reverse-H\"older inequality for Sobolev eigenfunctions,
in particular for the exponents $p-1$ and $p$, is a key step in our technique. 
We set out in Section~\ref{sec:revHol} the current state of play for 
reverse-H\"older inequalities in this context. 
It is tempting to ask for similar results in the remaining cases, 
when $p \neq r$, but we lack a reverse-H\"older inequality 
similar to \eqref{ineq:chiti2}. 

\bigskip \noindent {\sc Acknowledgements:} M.\ M.\ F.\ is partially supported by the 
Alexander vonHumboldt Foundation, and J.\ R.\ is partially supported by the National 
Research Foundation 
of South Africa. Part of this research was completed while M.\ M.\ F.\ visited J.\ R.\ at the 
University of Cape Town, and part of it while J.\ R.\ visited T.\ C.\ at University 
College Cork. We thank these institutions for their hospitality. 

\section {Hadamard variation formula}\label{sec:HadVar}
Following Grinfeld's approach in \cite[Section~5]{G}, 
we derive the Hadamard variation formula in a slightly more general 
setting than we require here. 
Take $X: (-\epsilon, \epsilon) \times \bar \Omega \rightarrow \R^n$ 
to be a time-dependent vector field on the closure of $\Omega$ 
and let $\xi: (-\epsilon, \epsilon) \times \bar \Omega \rightarrow \R^n$ 
be its flow, so that 
\[
\xi(0,x) = x ,\qquad \frac{\partial \xi}{\partial t} (t,x) = X(t,x). 
\]
Set $\Omega_t = \xi(t,\cdot )(\Omega)$. 

\begin {lemma} We have 
\begin {equation} \label {hadamard-var2} 
\dot {\mathcal{C}}_{p,r} 
	= \left. \frac{d}{dt} \mathcal{C}_{p,r}(\Omega_t)\right |_{t=0} 
	=(1-r) \int_{\partial \Omega} |\nabla \phi |^r 
		\langle X, \eta \rangle d\sigma. 
\end {equation} 
\end {lemma} 

\begin {proof}
For each $t \in (-\epsilon, \epsilon)$ we let $\phi_t$ be the extremal 
Sobolev function on $\Omega_t$, normalized 
so that $\int_{\Omega_t} \phi_t^p \,d\mu = 1$. 
Differentiating the normalization with respect to $t$ and using the fact that 
$\phi$ vanishes on the boundary of $\Omega$ gives 
\[
0 	= \left. \frac{d}{dt} \int_{\Omega_t} \phi_t^p \,d\mu \right |_{t=0} 
	= p \int_{\Omega} \phi^{p-1} \frac{\partial \phi}{\partial t} \,d\mu 
		+ \int_{\partial \Omega} \phi^p \langle X, \eta \rangle \,d\sigma 
	= p \int_{\Omega} \phi^{p-1} \frac{\partial \phi}{\partial t} \,d\mu
\]
so that 
\begin {equation} \label{variation1} 
\int_{\Omega} \phi^{p-1} \frac{\partial \phi}{\partial t} \,d\mu = 0.
\end {equation} 
Next we differentiate the boundary condition 
$\left. \phi_t \right |_{\partial \Omega_t} = 0$ with respect to $t$ at $t=0$ 
to obtain that 
\[
0 	= \frac{\partial \phi}{\partial t} + \langle X, \nabla \phi \rangle 
	=\frac{\partial \phi}{\partial t} - |\nabla \phi | \langle X,\eta\rangle,
\]
so that 
\begin {equation} \label{variation3} 
\frac{\partial \phi}{\partial t} =  |\nabla \phi| \langle X, \eta \rangle 
\qquad \textrm{on } \partial \Omega.
\end {equation} 
Thus 
\begin {align*} 
\dot {\mathcal{C}}_{p,r} 
& =  \frac{d}{dt} \int_\Omega 
	\big[ \langle \nabla \phi, \nabla \phi \rangle \big]^{r/2} \,d\mu \\ 
& =  r \int_\Omega |\nabla \phi|^{r-2} \left \langle \nabla \phi, 
	\nabla \frac{\partial \phi}{\partial t} \right \rangle \,d\mu 
	+ \int_{\partial \Omega} |\nabla \phi|^r \langle X, \eta \rangle \,d\sigma \\ 
& =  - r \int_\Omega \frac{\partial \phi}{\partial t} \Delta_r \phi \,d\mu 
	+ r \int_{\partial \Omega} \frac{\partial \phi}{\partial t} |\nabla \phi|^{r-2} 	
	\frac{\partial \phi}{\partial \eta}\,d\sigma 
	+ \int_{\partial \Omega} |\nabla \phi|^r \langle X, \eta \rangle \,d\sigma \\ 
& =  r \,\mathcal{C}_{p,r}(\Omega) \int_\Omega \phi^{p-1} 
	\frac{\partial \phi}{\partial t} \,d\mu 
	- r \int_{\partial \Omega} \frac{\partial \phi}{\partial t} |\nabla \phi|^{r-1} 
	\,d\sigma + \int_{\partial \Omega}  |\nabla \phi|^r \langle X, \eta \rangle 
	\,d\sigma \\ 
& =  -r \int_{\partial \Omega} |\nabla \phi|^r \langle X, \eta \rangle 
	\,d\sigma + \int_{\partial \Omega} |\nabla \phi|^r \langle X, \eta \rangle 
	\,d\sigma . \qedhere
\end {align*} 
\end {proof}

\section {Rearrangements}\label{sec:rearrangements}
We derive some preliminary rearrangement inequalities needed 
in Section~\ref{sec:revHol} to prove 
reverse-H\"older inequalities for the eigenfunctions $\phi$. 
We set $M = \sup _{x \in \Omega} (\phi(x))$ and, for $t \in [0,M]$, set 
\[
\Omega_t = \{ x \in \Omega: \phi(x) > t\}, \qquad \mu (t) = |\Omega_t|.
\]
This distribution function $\mu$ is nonincreasing, so it has an inverse 
function 
$$\phi^* : [0,|\Omega|] \rightarrow [0,M], \qquad \phi^*(v) = \inf \{ t \in 
[0,M]: \mu(t)>v\}.$$
Observe that both $\mu$ and $\phi^*$ are differentiable almost everywhere and 
(when they are both defined) we have 
$$\mu'(t) = \frac{1}{(\phi^*)'(\mu(t))} = - \int_{\{\phi = t\} }  \frac{d\sigma}{|\nabla \phi |}.$$

In the next section we will compare $\phi$ and $\phi^*$ to the corresponding 
extremal functions $\psi$ and $\psi^*$ on $\B^*$, the round ball with 
$\mathcal{C}_{p,r}(\Omega) = \mathcal{C}_{p,r}(\B^*)$, so we take 
this opportunity to record the equations which $\psi$ and $\psi^*$ satisfy. 
The function $\psi$ is radial and decreasing, so (see the introduction of \cite{Kaw})
\begin {eqnarray*} 
- \mathcal{C}_{p,r} (\Omega) \psi^{p-1} & = &  -\mathcal{C}_{p,r}(\B^*)  
= \Delta_r \psi \\ 
& = &  \left ( - \frac{d \psi}{d \rho} \right )^{r-2} \left [ 
(r-1) \frac{d^2 \psi}{d \rho^2} + \frac{n-1}{\rho} \frac{d \psi}
{d \rho} \right ] \\ 
& =&  - \rho^{1-n} \frac{d }{d \rho} \left ( 
\rho^{n-1} \left ( -\frac{d \psi}{ d\rho} \right )^{r-1} \right) . 
\end {eqnarray*} 
We change variables to $v = \omega_n \rho^n$, and define $\psi^*(v) 
= \psi\left ( \left ( \frac{v}{\omega_n} \right )^{1/n} \right )$. Then 
$$ \mathcal{C}_{p,r} (\psi^*(v))^{p-1} = n^r \omega_n^{r/n}  
\frac{d}{dv} \left [ v^{\frac{r(n-1)}{n}} \left (- \frac{d\psi^*}{dv} \right )^{r-1}
\right ],$$
which we can integrate once and rearrange to read 
\begin {equation} \label {ball-eigenfunction} 
\left ( - (\psi^*)'(v)\right )^{r-1} = n^{-r} \omega_n^{-r/n} 
v^{\frac{r(1-n)}{n}} \mathcal{C}_{p,r} (\Omega)\int_0^v (\psi^*(\tau))^{p-1} 
d\tau. \end {equation} 

The following is an adaptation of Talenti's inequality (see (34) of \cite{Tal}). 
\begin {lemma} 
We have 
\begin {equation} \label{talenti-ineq} 
(- (\phi^*)'(v))^{r-1} \leq n^{-r} \omega_n^{-r/n} \mathcal{C}_{p,r} 
(\Omega) v^{\frac{r(1-n)}{n}} \int_0^v 
(\phi^*(\tau))^{p-1} d\tau \end {equation}
for almost every $v$, where $\omega_n$ is the volume of a unit ball in $\R^n$.  
Moreover, equality can only occur if $\Omega$ is a round ball. 
\end {lemma} 

\begin {proof} The fact that $\phi$ is an extremal function implies 
$$\mathcal{C}_{p,r} (\Omega) \int_{\Omega_t} \phi^{p-1} d\mu = - \int_{\Omega_t} 
\Delta_r \phi d\mu = - \int_{\partial \Omega_t} |\nabla \phi|^{r-2} \frac{\partial \phi}
{\partial \eta} d\sigma = \int_{\partial \Omega_t} |\nabla \phi|^{r-1} d\sigma. $$
We combine this inequality with H\"older's inequality to obtain 
\begin {eqnarray*} 
|\partial \Omega_t| & = & \int_{\partial \Omega_t} d\sigma 
= \int_{\partial \Omega_t} |\nabla \phi |^{\frac{r-1}{r}} |\nabla \phi |^{\frac{1-r}{r}} 
d\sigma \\ 
& \leq & \left ( \int_{\partial \Omega_t} |\nabla \phi |^{r-1} d\sigma \right )^{1/r}
\left ( \int_{\partial \Omega_t} |\nabla \phi|^{-1} d\sigma \right )^{\frac{r-1}{r}} \\ 
& = & (-\mu'(t))^{\frac{r-1}{r}} \left ( \mathcal{C}_{p,r} (\Omega) \int_{\Omega_t} 
\phi^{p-1} d\mu \right )^{1/r} ,\end {eqnarray*} 
which we can rearrange to read  
$$\mathcal{C}_{p,r}(\Omega)(-\mu'(t))^{r-1}  \int_{\Omega_t} \phi^{p-1} 
d\mu \geq |\partial \Omega_t|^r \geq \left [ n \omega_n \left (
\frac{1}{\omega_n} \mu(t)  \right )^{\frac{n-1}{n}} \right ]^r.$$
The inequality \eqref{talenti-ineq} now follows once we change variables to $v = \mu(t)$ 
and recall   
$\mu' = \frac{1}{(\phi^*)'}$. Moreover, equality in \eqref {talenti-ineq} 
forces equality in our use of the isoperimetric inequality, which forces $\Omega$ 
to be a round ball. 
\end {proof} 

It is crucial that the right hand sides of \eqref{ball-eigenfunction} and 
\eqref{talenti-ineq} are essentially the same. 

\section{Reverse-H\"older inequalities}\label{sec:revHol}

In this section we prove inequalities which reverse the usual H\"older inequality 
for extremal functions $\phi$ in several cases. We summarize our results with the following
theorem. 
\begin {thm} Let $\phi$ be an extremal function and let $0<q_1<q_2$. A 
reverse-H\"older inequality of the form 
$$\left ( \int_\Omega \phi^{q_1} d\mu \right )^{1/q_1} \geq C \left ( 
\int_\Omega \phi^{q_2} d\mu \right )^{1/q_2}$$
holds in the following cases: 
\begin {itemize} 
\item $p=r$,
\item $n=2$ and $r=2$, $q_1=p-1$ and $q_2 = p$,
\item $q_1 = p$.
\end {itemize} 
In all cases the constant $C$ depends only on $n,p,r,q_1, q_2$. Moreover  
equality implies $\Omega$ is a round ball. \end {thm} 
We prove the case of $p=r$ below in Proposition \ref{prop:chiti2}, 
and prove the case of $q_1 = p$ in Proposition \ref{prop:chiti3}. One 
can find a proof of the case $n=2, r=2, q_1 = p-1, q_2 = p$ in 
\cite{CR1}. 

Below we will see that the proof of Proposition \ref{prop:chiti2} is 
easier than the proof of Proposition \ref{prop:chiti3}, mostly because 
\eqref{sobolev-pde} is homogeneous only in the case $p=r$. The homogeneity 
allows us to multiply $\psi$ and 
$\phi$ by convenient constants, so that we can choose a scale on which to work. 
It is curious to us that in this particular application homogeneity is even more 
important than linearity. 

In most of our computations for thise section we will temporarily drop the 
normalizations $$\int_\Omega \phi^p d\mu = 1, \qquad 
\int_{\B^*} \psi^p d\mu = 1.$$

We will compare $\Omega$ to $\B^*$, the round ball with 
$\mathcal{C}_{p,p}(\Omega) = \mathcal{C}_{p,p}(\B^*)$. An important tool we 
use is the Faber-Krahn inequality, which implies $|\Omega| \geq |\B^*|$, with 
equality if and only if $\Omega = \B^*$. 

\begin {prop} \label{prop:chiti1}
Let $\phi$ be an extremal function on $\Omega$ and let $\psi$ 
be the extremal function of $\B^*$, the ball with $\mathcal{C}_{p,p} (\B^*) 
= \mathcal{C}_{p,p}(\Omega)$.  Normalized both $\phi$ and $\psi$ 
so that $\| \phi \|_{L^\infty}
= \| \psi \|_{L^\infty}$. Then for $0 < v < |\B^*|$ we have $\phi^*(v) \geq \psi^*(v)$. 
Moreover, equality can occur for some $v>0$ only if $\Omega = \B^*$. \end {prop} 

\begin {proof} If $|\Omega| = |\B^*|$ then $\Omega = \B^*$ by the Faber-Krahn 
inequality, and in this case there is nothing to prove, so we assume $|\Omega| 
> |\B^*|$. In this case 
$$\phi^* (0) = \psi^*(0) = \| \phi\|_{L^\infty}, \qquad \psi^*(|\B^*|) = 0 < \phi^*(|\B^*|),$$
so there must exist $k>1$ such that $k\phi^*(v) \geq \psi^*(v)$ on the interval 
$[0,|\B^*|]$. We define 
$$k_0 = \inf \{ k > 1 : k \phi^*  \geq \psi^* \}$$ 
and complete the proof by showing $k_0 = 1$. If $k_0 > 1$ then 
there exists $v_0 \in (0,|\B^*|)$ such that $\psi^*(v_0) = k_0\phi^*(v_0)$ and 
$\psi^*(v) < k_0\phi^*(v)$ on $(0,v_0)$. Now define 
$$u^*: [0, |\B^*|] \rightarrow \R, \qquad u^*(v) = \left \{ \begin {array}{rl} 
k_0 \phi^*(v) & 0 \leq v \leq v_0 \\ \psi^*(v) & v_0 \leq v \leq |\B^*|  \end {array} 
\right. $$
and let $u(x) = u^*( \omega_n |x|^n)$. Then by \eqref{ball-eigenfunction} 
and \eqref{talenti-ineq} we have 
$$(-(u^*)'(v))^{p-1}  \leq n^{-r} \omega_n^{-p/n} v^{\frac{p(1-n)}{n}} \mathcal{C}_{p,p}
(\Omega) \int_0^v (u^*(\tau))^{p-1} d\tau,$$
and so 
\begin {eqnarray*} 
\int_{\B^*} |\nabla u|^p d\mu & = & \int_{\B^*} \left ( \frac{du}{d\rho} \right )^p d\mu 
= \int_0^{|\B^*|} n^p \omega_n^{p/n} v^{\frac{p(n-1)}{n}} ((u^*)'(v))^p dv \\ 
& \leq & \mathcal{C}_{p,p}(\Omega)  \int_0^{|\B^*|} (-(u^*)'(v))
\int_0^v (u^*(\tau))^{p-1} d\tau dv \\ 
& = & \mathcal{C}_{p,p}(\Omega) \int_0^{|\B^*|}(u^*(\tau))^{p-1} \int_\tau^{|\B^*|} (-(u^*)'(v)) dv d\tau\\ 
& = & \mathcal{C}_{p,p}(\Omega) \int_{\B^*} u^p d\mu . \end {eqnarray*} 
However, $\mathcal{C}_{p,p}(\B^*) = \mathcal{C}_{p,p}(\Omega)$, so 
this is only possible if $u = \psi$, which cannot occur because $u^* > \psi^*$ 
on $(0,v_0)$. 
\end {proof}

\begin {cor} For any $q \geq 0$ we have the scale-invariant inequality 
\begin {equation} \label{ineq:chiti1}
\frac{\|  \phi\|_{L^q}}{\| \phi \|_{L^\infty}} \geq 
\frac{ \| \psi \|_{L^q}} {\| \psi \|_{L^\infty}}\end {equation} 
with equality if and only if $\Omega = \B^*$. \end {cor} 
\begin {proof} Integrate the inequality we've just proved in Proposition \ref{prop:chiti1}.
\end {proof} 

\begin {prop} \label{prop:chiti2}
Let $0 < q_1 < q_2 < \infty$. There exists $K$ depending only on 
$n$, $p$, $q_1$, and $q_2$ such that 
$$\left ( \int_\Omega \phi^{q_1} d\mu \right )^{q_2} \geq K 
(\mathcal{C}_{p,p}(\Omega))^{-\frac{n}{p} (q_2 - q_1)}
\left ( \int_\Omega \phi^{q_2} d\mu \right )^{q_1},$$
and equality implies $\Omega$ is a round ball. 
\end {prop} 

\begin {proof} If $|\Omega| = |\B^*|$ there is nothing to prove, so 
we assume $|\Omega| > |\B^*|$. This time we choose the normalization 
$$\int_\Omega \phi^{q_1} d\mu = \int_{\B^*} \psi^{q_1} d\mu,$$
so that \eqref{ineq:chiti1} implies 
$$ \psi^*(0) = \| \psi\|_{L^\infty} > \| \phi \|_{L^\infty} = \phi^*(0).$$
We also know that $\phi^*(|\B^*|) > 0 = \psi^*(|\B^*|)$, so the graphs of 
the functions $\phi^*$ and $\psi^*$ must cross somewhere in the interval 
$(0, |\B^*|)$. Let 
$$v_0 = \inf \{ v \in [0, |\B^*|]: \phi^*(\tilde v) \leq \psi^*(\tilde v) 
\textrm{ for all } \tilde v \in (0,v) \}$$
be the first crossing when viewed from the left. Then 
$$0 < v_0 < |\B^*|, \quad \psi^*\geq \phi^* \textrm{ in }[0,v_0], \quad 
\psi^*(v_0) = \phi^*(v_0)$$
and there exists $v  \in (v_0, |\B^*|)$ such that $\phi^*(v) > \psi^*(v)$. In 
fact, by continuity the inequality $\phi^*> \psi^*$ must hold in a nontrivial 
interval $I$ surrounding $v$.  

We claim that $\phi^* > \psi^*$ on the entire interval $(v_0, |\B^*|)$. Suppose 
otherwise, then there must exist $v_1 \in (v_0, |\B^*|)$ with $\phi^*(v_1) = 
\psi^*(v_1)$ and we can define 
$$u^*(v) = \left \{ \begin {array}{rl} \psi^*(v) & 0 \leq v \leq v_0 \\ 
\max \{ \psi^*(v), \phi^*(v) \} & v_0 \leq v \leq v_1 \\ \psi^*(v) & v_1 \leq v 
\leq |\B^*| \end {array} \right. $$ 
and $u(x) = u^*(\omega_n |x|^n)$. 
Again by \eqref{ball-eigenfunction} and \eqref{talenti-ineq} we have 
$$(-(u^*)'(v))^{p-1}  \leq n^{-p} \omega_n^{-p/n} v^{\frac{p(1-n)}{n}} \mathcal{C}_{p,p}
(\Omega) \int_0^v (u^*(\tau))^{p-1} d\tau,$$
so, as in our proof of Proposition \ref{prop:chiti1} we have 
\begin {eqnarray*}
\int_{\B^*} |\nabla u|^p d\mu & = & \int_{\B^*} \left ( \frac{du}{d\rho} \right )^p d\mu 
= \int_0^{|\B^*|} n^p \omega_n^{p/n} v^{\frac{p(n-1)}{n}} ((u^*)'(v))^p dv \\ 
& \leq & \mathcal{C}_{p,p}(\Omega)  \int_0^{|\B^*|} (-(u^*)'(v))
\int_0^v (u^*(\tau))^{p-1} d\tau dv \\ 
& = & \mathcal{C}_{p,p}(\Omega) \int_0^{|\B^*|}(u^*(\tau))^{p-1} \int_\tau^{|\B^*|} (-(u^*)'(v)) dv d\tau\\ 
& = & \mathcal{C}_{p,p}(\Omega) \int_{\B^*} u^p d\mu . \end {eqnarray*} 
That $\mathcal{C}_{p,p}(\Omega) = \mathcal{C}_{p,p}(\B^*)$ now implies 
$u$ is a muliple of $\psi$, which is impossible. 

We conclude that $\psi^* \geq \phi^*$ on $[0, v_0]$ and $\psi^* < \phi^*$ 
on $(v_0, |\B^*|]$. Then the argument in Theorem 7 of \cite{CR2} shows 
$$\left ( \int_\Omega \phi^{q_2} d\mu \right )^{1/q_2} \leq \left ( \int_{\B^*} 
\psi^{q_2} d\mu \right )^{1/q_2} = \frac{\left ( \int_{\B^*} \psi^{q_2} d\mu \right )^{1/q_2}}
{ \left ( \int_{\B^*} \psi^{q_1} d\mu \right )^{1/q_1}} \left ( \int_\Omega \phi^{q_1} d\mu 
\right )^{1/q_1},$$
which we can rewrite as 
$$ \left ( \int_\Omega \phi^{q_1} d\mu \right )^{q_2} \geq \widetilde C \left (
\int_\Omega \phi^{q_2} d\mu \right )^{q_1}, \qquad \widetilde C = \frac
{\left ( \int_{\B^*} \psi^{q_1} d\mu \right )^{q_2}}{\left (\int_{\B^*} \psi^{q_2}
d\mu \right )^{q_1}}.$$

All that remains now is to unravel the constant $\widetilde C$. Let $R$ be the 
radius of $\B^*$ and define the function 
$$\widetilde \psi : \B_1 \rightarrow \R, \qquad \widetilde \psi(x) = \psi(Rx).$$
Then $\widetilde \psi$ is an extremal function for $\mathcal{C}_{p,p}(\B_1)$, because 
the PDE \eqref{sobolev-pde} is homogeneous in the case $p=r$. 
By \eqref{sobolev-scaling} we have 
$$\mathcal{C}_{p,p}(\B^*) = \mathcal{C}_{p,p} (\B_R) = R^{-p} \mathcal{C}_{p,p} (\B_1)$$
which implies 
$$R = \left ( \frac{\mathcal{C}_{p,p} (\B^*)}{\mathcal{C}_{p,p}(\B_1)} 
\right )^{-1/p} = \left ( \frac{\mathcal{C}_{p,p} (\Omega)}{\mathcal{C}_{p,p}(\B_1)} 
\right )^{-1/p} ,$$ 
and so 
$$\widetilde C = R^{n(q_2 - q_1)} \frac{\left ( \int_{\B_1} \widetilde \psi^{q_1} d\mu 
\right )^{q_2}}{\left ( \int_{\B_1} \widetilde \psi^{q_2} d\mu \right )^{q_1}}
= \left ( \frac{\mathcal{C}_{p,p} (\Omega)}{\mathcal{C}_{p,p}(\B_1)} 
\right )^{-\frac{n}{p} (q_2 - q_1)} \frac{\left ( \int_{\B_1} \widetilde \psi^{q_1} d\mu 
\right )^{q_2}}{\left ( \int_{\B_1} \widetilde \psi^{q_2} d\mu \right )^{q_1}}.$$
\end {proof} 

We will use the case of $q_1 = p-1$ and $q_2 = p$ in the next section. 
\begin {cor} There exists a constant $K$ depending only on $n$ and $p$ such that 
\begin {equation} \label {ineq:chiti2} 
\left ( \int_\Omega \phi^{p-1} d\mu \right )^p \geq K (\mathcal{C}_{p,p}
(\Omega))^{-n/p} \left ( \int_\Omega \phi^p d\mu \right )^{p-1}. \end {equation} 
Equality can only occur if $\Omega$ is a round ball. 
\end {cor} 

We close this section with a result generalizing the main theorem of \cite{CR3}. 
\begin {prop}  \label{prop:chiti3} 
Let $1 \leq r <n$, $1 \leq p < \frac{nr}{n-r}$, and $q>p$. There 
exists $K>0$ depending only on $n,r,p,q$ such that 
\begin {equation}\label{ineq:chiti3}
\left ( \int_\Omega \phi^p d\mu \right )^{1/p} \geq K 
(\mathcal{C}_{p,r}(\Omega))^{\frac{n(q-p)}{p(np-rp-nr)}}
\left (\int_\Omega \phi^q d\mu \right )^{1/q} \end {equation} 
for all extremal functions $\phi$. 
\end {prop} 

\begin {proof} As before we let $\B^*$ be the ball with $\mathcal{C}_{p,r} (\Omega) 
= \mathcal{C}_{p,r}(\B^*)$, and let $\psi$ be the corresponding extremal function 
on $\B^*$. By the Faber-Krahn inequality, we have $|\Omega |\geq |\B^*|$, with equality 
if and only if $\Omega = \B^*$. In the case $|\Omega| = |\B^*|$ we must also have 
equality in \eqref{ineq:chiti3}, which will more precisely read 
$$\left ( \int_\Omega \phi^p d\mu \right )^{1/p} \geq \widetilde K \left ( \int_\Omega 
\phi^q d\mu \right )^{1/q}, \qquad \widetilde K = \frac{\left ( \int_{\B^*} \psi^p d\mu 
\right )^{1/p}}{\left ( \int_\Omega \psi^q d\mu \right )^{1/q}};$$
we will see that in fact this constant $\widetilde K$ is optimal in general. 
Furthermore, if we let $R$ be the radius of $\B^*$ then \eqref{sobolev-scaling} 
implies 
$$R = \left ( \frac{\mathcal{C}_{p,r}(\B^*)}{\mathcal{C}_{p,r}(\B)} 
\right )^{\frac{p}{np-rn-rp}},$$
so that 
\begin {eqnarray*}
\widetilde K & = & \frac{\left ( \int_{\B^*} \psi^p d\mu 
\right )^{1/p}}{\left ( \int_{\B^*} \psi^q d\mu \right )^{1/q}} = R^{\frac{n(q-p)}{qp}}
\frac{\left ( \int_{\B} \bar \psi^p d\mu 
\right )^{1/p}}{\left ( \int_\B \bar \psi^q d\mu \right )^{1/q}} \\ 
& = & \left ( \frac{\mathcal{C}_{p,r}(\B^*)}{\mathcal{C}_{p,r}(\B)} \right )^{\frac{n(q-p)}
{q(np-rp-nr)}}\frac{\left ( \int_{\B} \bar \psi^p d\mu 
\right )^{1/p}}{\left ( \int_\B \bar \psi^q d\mu \right )^{1/q}} =
K \left ( \mathcal{C}_{p,r}(\Omega) \right )^{\frac{n(q-p)}{q(np-rp-nr)}} ,
\end {eqnarray*}
where $\bar \psi: \B \rightarrow \R$, $\bar \psi(x) = \psi (Rx)$ is the 
extremal function on the unit ball $\B$. 

Next we treat the case $|\Omega | > |\B^*|$. Normalize both extremal functions $\phi$ and 
$\psi$ so that 
$$\int_\Omega \phi^p d\mu = 1, \qquad \int_{\B^*} \psi^p d\mu = 1.$$ 
Combining these normalizations with $|\Omega| > |\B^*|$ we see 
\begin {equation} \label{crossing1}
1 = \int_0^{|\B^*|} (\psi^*)^p dv = \int_0^{|\Omega|}(\phi^*)^p dv > \int_0^{|\B^*|}
(\phi^*)^p dv,\end {equation} 
which implies we cannot have $\psi^*\leq\phi^*$ on the whole of the interval $[0,|\B^*|]$.  
On the other hand, we know 
$$\psi^*(|\B^*|) = 0 < \phi^*(|\B^*|),$$
so the graphs of these two functions must cross. Define 
$$v_1 = \inf\{ v \in [0,|\B^*|]: \psi^*(\tilde v) < \phi^*(\tilde v)
\textrm{ for all }\tilde v \in (v,|\B^*|] \}; $$
this is the first crossing of the two graphs, when viewed from the 
right hand side. By continuity, $\psi^*(v_1) = \phi^*(v_1)$ and $\psi^* 
< \phi^*$ on the interval $(v_1, |\B^*|]$. We also cannot have $v_1=0$, 
as this would contradict \eqref{crossing1}. 

We claim that $\psi^* \geq \phi^*$ on the interval $[0,v_1]$. Indeed, if this 
inequality did not hold, then we must have $\psi^*(v_2) < \phi^*(v_2)$ for 
some $v_2 \in [0,v_1)$, and by continuity this inequality must extend to 
an interval containing $v_2$. 

The function 
$$w^*:[0,|\B^*|] \rightarrow [0,\infty), \qquad w^*(v) = \left \{ \begin {array}{rl}
\max\{\phi^*(v), \psi^*(v)\} & 0\leq v \leq v_1\\ \psi^*(v) & v_1\leq v \leq |\B^*|\\
\end {array} \right. $$
satisfies 
$$(-(w^*)'(v) )^{r-1} \leq n^{-r} \omega_n^{-r/n} \mathcal{C}_{p,r}(\Omega)
v^{\frac{r(1-n)}{n}} \int_0^v (w^*(\tau))^{p-1} d\tau $$
by \eqref{ball-eigenfunction} and \eqref{talenti-ineq}. Now we can 
define 
$$w:\B^* \rightarrow \R, \qquad w(x) = w^*(\omega_n |x|^n),$$
so that we have 
\begin {eqnarray*} 
\int_{\B^*} |\nabla w|^r d\mu & = & \int_{\B^*} \left ( \frac{dw}{d\rho} \right )^r d\mu 
= \int_0^{|\B^*|} n^r \omega_n^{r/n} v^{\frac{r(n-1)}{n}} ((w^*)'(v))^r dv \\ 
& \leq & \mathcal{C}_{p,r}(\Omega)  \int_0^{|\B^*|} (-(w^*)'(v))
\int_0^v (w^*(\tau))^{p-1} d\tau dv \\ 
& = & \mathcal{C}_{p,r}(\Omega) \int_0^{|\B^*|} (w^*(\tau))^{p-1} \int_\tau^{|\B^*|} (-(w^*)'(v)) dv d\tau\\ 
& = & \mathcal{C}_{p,r}(\Omega) \int_{\B^*} w^p d\mu . \end {eqnarray*} 
Since $\mathcal{C}_{p,p}(\Omega) = \mathcal{C}_{p,p}(\B^*)$, the function $w$ must 
be extremal on $\B^*$, which implies $w$ must be a scalar mulitple of $\psi$. This 
would contradict the fact that $\psi<
\phi$ on an interval containing $v_2$. 

We've concluded that the graphs of $\phi^*$ and $\psi^*$ cross exactly once 
on the interval $[0,|\B^*|]$. The remainder of the argument is exactly the 
same as the one in \cite{CR3}, and we refer the reader to this treatment. 
\end {proof}

\section {Proof of Theorem \ref{thm:n=r=2}} 

The proof of Theorem \ref{thm:n=r=2} is easiest, so we present it first. 

We will need a reverse-H\"older inequality proved in \cite{CR1}, 
which reads 
\begin {equation} \label {reverse-holder}
\left ( \int_\Omega \phi^{p-1} d\mu \right )^2 \geq \frac{8\pi}{p 
\mathcal{C}_{p,2}(\Omega)} \left ( \int_\Omega \phi^p d\mu \right 
)^{\frac{2p-2}{p}}.\end {equation} 

In our setting $\left. X \right |_{\partial \Omega} = e^w \eta$, so 
\begin {eqnarray*} 
-\dot {\mathcal{C}}_{p,2} & = & \int_{\partial \Omega} e^w \left ( 
\frac{\partial \phi}{\partial \eta} \right )^2 d\sigma \geq \frac{1}
{\int_{\partial \Omega} e^{-w} d\sigma} \left ( \int_{\partial \Omega} 
\frac{\partial \phi}{\partial \eta}  d\sigma \right )^2 \\ 
& = & \frac{1}{\int_{\partial \Omega} e^{-w} d\sigma} 
\left ( \int_\Omega \Delta \phi d\mu \right )^2 = \frac{(\mathcal{C}_p)^2}
{\int_{\partial \Omega} e^{-w} d\sigma} \left ( \int_\Omega \phi^{p-1} d\mu \right )^2 \\ 
& \geq & \frac{8\pi \mathcal{C}_p}{p \int_{\partial \Omega} e^{-w} d\sigma}
\left ( \int_\Omega \phi^p d\mu \right )^{\frac{2(p-1)}{p}} = \frac{8\pi 
\mathcal{C}_p}{p \int_{\partial \Omega} e^{-w} d\sigma},\end {eqnarray*}
which proves \eqref{ineq:n=r=2}. Here we have first used \eqref{hadamard-var2}, 
followed by the Cauchy-Schwarz inequality, the divergence theorem, \eqref{sobolev-pde}, 
\eqref{reverse-holder}, and \eqref{sobolev-normalization}. Furthermore, equality 
in \eqref{ineq:n=r=2} forces equality in all the inequalities we have used. Equality 
in our use of \eqref{reverse-holder} can only occur if $\Omega$ is a round disk 
and equality in our use of the Cauchy-Schwarz inequality can only occur if $w$ 
is constant. \hfill $\square$

\section {Proof of Theorem \ref{thm:p=r}}

We first observe that 
$$|\nabla \phi|^{p-1} = e^{\frac{w(p-1)}{p}} |\nabla \phi|^{p-1} \cdot 
e^{\frac{w(1-p)}{p}},$$
so H\"older's inequality with exponents $\frac{p}{p-1}$ and $p$ 
gives us 
$$\int_{\partial \Omega} |\nabla \phi|^{p-1} d\sigma \leq 
\left ( \int_{\partial \Omega} e^w |\nabla \phi|^p d\sigma \right )^{\frac{p-1}{p}}
\left ( \int_{\partial \Omega} e^{(1-p)w} d\sigma \right )^{1/p}, $$
which we can rewrite as 
$$
\int_{\partial \Omega} e^w |\nabla \phi|^p d\sigma 
\geq \frac{\left ( \int_{\partial \Omega} |\nabla \phi|^{p-1} d\sigma 
\right )^{\frac{p}{p-1}} }{\left ( \int_{\partial \Omega} e^{(1-p)w} d\sigma 
\right )^{p-1} }. $$
Thus \eqref{hadamard-var2}, H\"older's inequality, the divergence theorem, and 
\eqref{ineq:chiti2} combine to give us 
\begin {eqnarray} \label{variation1}
-\dot {\mathcal{C}}_{p,p} & = & (p-1) \int_{\partial \Omega} e^w |\nabla \phi|^p 
d\sigma \\ \nonumber 
& \geq & \frac{p-1}{\left ( \int_{\partial \Omega} e^{(1-p) w}d\sigma \right )^{p-1}}
\left ( \int_{\partial \Omega} |\nabla \phi |^{p-1} d\sigma \right )^{\frac{p}{p-1}} \\ \nonumber 
& = & \frac{p-1}{\left ( \int_{\partial \Omega} e^{(1-p) w}d\sigma \right )^{p-1}}
\left ( \int_{\partial \Omega} -|\nabla \phi |^{p-2} \frac{\partial \phi}{\partial \eta} 
d\sigma \right )^{\frac{p}{p-1}} \\ \nonumber 
& = & \frac{p-1}{\left ( \int_{\partial \Omega} e^{(1-p) w}d\sigma \right )^{p-1}}
\left ( \int_\Omega -\Delta_p \phi  d\mu \right )^{\frac{p}{p-1}} \\ \nonumber 
& = & \frac{(p-1) (\mathcal{C}_{p,p})^{\frac{p}{p-1}}}
{\left ( \int_{\partial \Omega} e^{(1-p) w}d\sigma \right )^{p-1}}
\left [ \left ( \int_\Omega  \phi^{p-1}  d\mu \right )^p \right ]^{\frac{1}{p-1}}\\  \nonumber 
& \geq & \frac{(p-1) (\mathcal{C}_{p,p})^{\frac{p}{p-1}}}
{\left ( \int_{\partial \Omega} e^{(1-p)w}d\sigma\right )^{p-1}} 
\left [ K (\mathcal{C}_{p,p})^{-n/p}  \left ( \int_\Omega
\phi^p d\mu \right )^{p-1} \right ]^{\frac{1}{p-1}} \\ \nonumber 
& = & \frac{(p-1)K^{\frac{1}{p-1}}}{\left ( \int_{\partial \Omega} 
e^{(1-p)w} d\sigma \right )^{p-1}} (\mathcal{C}_{p,p}(\Omega))^{\frac{1}{p-1} \left ( 
p - \frac{n}{p} \right )} . 
\end {eqnarray} 
In the case $1<p<n$ we can rearrange \eqref{variation1} to read 
$$ (\mathcal{C}_{p,p} )^{\frac{1}{p-1} \left ( \frac{n}{p} - p \right )} \dot {\mathcal{C}}_{p,p}
\geq \frac{(p-1) K^{\frac{1}{p-1}}} {\left ( \int_{\partial \Omega} e^{(1-p)w} d\sigma
\right )^{p-1}},$$
which implies \eqref{ineq:p=r}. In the case $p=n$ we rewrite \eqref{variation1} 
as 
$$- \frac{\dot {\mathcal{C}}_{n,n}}{\mathcal{C}_{n,n}} \geq \frac{(n-1) K^{\frac{1}{n-1}}}
{\left ( \int_{\partial \Omega} e^{(1-n)w} d\sigma \right )^{n-1}}$$
which implies \eqref{ineq:p=r=n}. Moreover, equality in \eqref{variation1} 
forces equality in \eqref{ineq:chiti2}, 
which in turn forces $\Omega$ to be a round ball. Also, equality in our use of the 
H\"older inequality can only occur if $e^{\frac{w(p-1)}{p}}$ is a multiple of $e^{\frac{w(1-p)}{p}}$, 
which forces $w$ to be constant. 
\hfill $\square$

\section {Proof of Theorem \ref{thm:conformal}}

In this setting we let $F: \B \rightarrow \R^n$ be a conformal diffeomorphism, and 
for $0<t<1$ we consider compare the balls $\B_t$ to their conformal images 
$\Omega_t = F(\B_t)$. Next we define $\mathcal{C}_{p,p}(t) = \mathcal{C}_{p,p} 
(\B_t)$, with its associated extremal function $\phi_t$, and 
$\widetilde {\mathcal{C}}_{p,p}= \mathcal{C}_{p,p}(\Omega_t)$, with its associated 
extremal function $\tilde \phi_t$. As usual, we choose the normalization 
$$\int_{\B_t} \phi^p d\mu = 1 = \int_{\Omega_t} (\tilde \phi)^p d\tilde \mu.$$
We also set $\psi = \tilde \phi \circ F$ and notice that, because $F$ is 
conformal, $|\nabla \psi| = |DF| |\nabla \tilde \phi|$. 
In this setting \eqref{hadamard-var2} reads 
\begin {equation} \label{hadamard-var3} 
\frac{d}{dt} \widetilde {\mathcal{C}}_{p,p} = (1-p) \int_{\partial \Omega_t}
|DF| |\nabla \tilde \phi|^p d\tilde \sigma = (1-p) \int_{\partial \B_t} |DF|^{n-2} 
|\nabla \psi|^p d\sigma.\end {equation} 

Combining our normalization with \eqref{ineq:chiti2}, we see 
\begin {eqnarray*} 
K \widetilde {\mathcal{C}}_{p,p}^{-n/p} & = & K \widetilde {\mathcal{C}}_{p,p}^{-n/p}
\left ( \int_{\Omega_t} (\tilde \phi)^p d\tilde \mu \right )^{p-1} \\
& \leq & \left ( \int_{\Omega_t} (\tilde \phi)^{p-1} d\tilde \mu \right )^p = 
(\widetilde {\mathcal{C}}_{p,p})^{-p} \left ( - \int_{\Omega_t} \Delta_p \tilde \phi 
d \tilde \mu \right )^p \\ 
& = & (\widetilde {\mathcal{C}}_{p,p})^{-p} \left ( \int_{\partial \Omega_t}
|\nabla \tilde \phi|^{p-1} d\tilde \sigma \right )^p,
\end {eqnarray*} 
which we can rewrite as 
\begin {eqnarray} \label{conf-var1} 
K (\widetilde {\mathcal{C}}_{p,p})^{p - \frac{n}{p}} & \leq & 
\left ( \int_{\partial \Omega_t} |\nabla \tilde \phi|^{p-1} d\tilde \sigma 
\right )^p = \left ( \int_{\partial \B_t} |DF|^{n-2}|\nabla \psi|^{p-1} d\sigma 
\right )^p \\ \nonumber 
& = & \left ( \int_{\partial \B_t} |DF|^{\frac{n-2}{p}} |DF|^{\left ( \frac{p-1}{p}
\right )(n-2)} |\nabla \psi|^{p-1} d\sigma \right )^p \\ \nonumber 
& \leq & \left ( \int_{\partial \B_t} |DF|^{n-2} |\nabla \psi|^p d\sigma \right )^{p-1}
\int_{\partial \B_t} |DF|^{n-2} d\sigma \\ \nonumber 
& = & \left ( \frac{1}{1-p} \frac{d}{dt} \widetilde {\mathcal{C}}_{p,p} \right )^{p-1} 
\int_{\partial \B_t} |DF|^{n-2} d\sigma. \end {eqnarray}
If $1<p<n$ then \eqref{conf-var1}  reads 
$$- \frac{d}{dt} (\widetilde {\mathcal{C}}_{p,p})^{\frac{n-p}{p(p-1)}} 
\geq \frac{ \left ( \frac{n-p}{p} \right ) K^{\frac{1}{p-1}}} 
{\left ( \int_{\partial \B_t} |DF|^{n-2} d\sigma \right )^{\frac{1}{p-1}}},$$
and \eqref{ineq:conf-p=r} follows from the equality case of \eqref{ineq:p=r} and 
the inequality $$\int_{\partial \B_t} |DF|^{n-2} d\sigma \geq |\partial \B_t|^{(p-1)^2}.$$ 
The proof of \eqref{ineq:conf-p=r=n} is very similar. \hfill $\square$

\end{document}